\documentclass[12pt,twoside]{amsart}

\usepackage{amssymb,amsmath}
\theoremstyle{plain}
\newtheorem{theorem}{Theorem}[section]
\newtheorem{lemma}[theorem]{Lemma}

\newtheorem{corollary}[theorem]{Corollary}
\newtheorem{question}[theorem]{Question}

\theoremstyle{definition}
\newtheorem{remark}[theorem]{Remark}

\numberwithin{equation}{section}

  \input xypic
  \xyoption{all}

 \usepackage{tikz}
\usetikzlibrary{arrows}

 \usepackage{times}

  \usepackage{amsmath, amsthm, amsfonts}

\usepackage{amssymb}
\usepackage{amsmath}
\usepackage{epsfig}
\usepackage{amsthm}
\usepackage{euscript}
\usepackage{latexsym, mathrsfs}

\numberwithin{equation}{section}
\def \mff{\mathsf}

\def \df{\Delta'}
\def \mc{\mathcal}
\def \inv{^{-1}}

\def \v{\vskip 0.1in}
\def \n{\noindent}

\def \b{\bar }

\def \p{\partial}

\begin{document}

\title{A Liouville Theorem on the PDE $\det(f_{i\bar j})=1$ }
\author[Li]{An-Min Li}
\address{Department of Mathematics,
Sichuan University,
 Chengdu, 610064, China}
\email{anminliscu@126.com}
\author[Sheng]{Li Sheng}
\address{Department of Mathematics,
Sichuan University,
Chengdu, 610064, China}
\email[Corresponding author]{lshengscu@gmail.com}
\thanks{ Li acknowledges the support of NSFC Grant NSFC11521061. \\
${}\quad\ $Sheng acknowledges the support of NSFC Grant NSFC11471225.}

\maketitle

\centerline{ Dedicated to Udo Simon   for his 80th Birthday}

{\abstract

 Let $f$  be a smooth plurisubharmonic function which  solves
􏰀$$ \det(f_{i\bar j})=1\;\;\;\;\;\;\mbox{in }\Omega\subset \mathbb C^n.$$ Suppose that the metric $\omega_{f}=\sqrt{-1}f_{i\bar j}dz_{i}\wedge d\bar z_{j}$ is complete and $f$ satisfies the growth condition
$$
C^{-1}(1+|z|^2)\leq f\leq C(1+ |z|^2),\;\;\;\; as\;\;\; |z|\to \infty.
$$
for some $C>0,$
then  $f$ is quadratic.
\endabstract
}



\section{\bf Introduction}\label{Sec-Intro}
Let $\Omega \subset \mathbb C^n$, denote
$$\mc R^\infty(\Omega):=\{f\in C^{\infty}(\Omega)\;|\; f\mbox{ is a
real function and }(f_{i \bar j})>0 \},
$$ where $(f_{i\bar j})=\left(\frac{\partial^2  f}
{\partial z_i\partial \bar{z}_j}\right).$ For $f\in \mc R^\infty(\Omega),$  $(\Omega, \omega_f)$ is
 a  K\"{a}hler  manifold. In this paper we study the PDE
\begin{equation}\label{eqn_1.1}
\det(f_{i\bar j})=1.
\end{equation}
 The quadratical polynomial
\begin{equation}\label{quadric}
f= z_1\bar z_1 + ...+ z_n\bar z_n\end{equation}
is a solution of \eqref{eqn_1.1}. If we take a linear transformation
 $$\dot{z}_i=\sum a_i^jz_j,\;\;\;\det(a_i^j)=1,$$
the above $f$ is transformed into another solution of \eqref{eqn_1.1}. All these solutions are called affine equivalent.

 When $\omega_{f}$ is complete, $\omega_{f}$ is a complete Calabi-Yau metric.
  Tian  proved that every Calabi-Yau metric of Euclidean volume growth on $\mathbb C^{2}$ must be Eucildean metric,  and conjectured that the same should hold true on $\mathbb C^n$ for all $n\geq 3$ (\cite{T}).
 Recently,   a counterexample to this conjecture was found by    Li, Conlon-Rochon,   Szekelyhidi independently (see \cite{CR,Li,Sz}).

 \v  To obtain the Liouville theorem of  \eqref{eqn_1.1}  people need to strengthen the assumption near $\infty.$
At the AIM workshop ``Nonlinear PDEs in real and complex geometry"  Szekelyhidi asked the following question (\cite{AIM})
\begin{question}
Let $f$ be plurisubharmonic solution to $det(f_{i\bar j})=1$ on $\mathbb C^n$ satisfying
$$C^{-1} (|z|^2+1)\leq f\leq C(|z|^2+1)$$
then $f$ is quadratic.
\end{question}
In \cite{W} Wang prove the Liouville theorem under the assumption $f=|z|^2+o(|z|^2),$ as $|z|\to \infty.$

\v
In this paper we consider the solution of \eqref{eqn_1.1} with complete K\"{a}hler metrics.
We prove the following two Liouville properties on
\eqref{eqn_1.1}.
\begin{theorem}\label{theorem_1.1} Let $f\in \mc R^\infty(\Omega)$ satisfying \eqref{eqn_1.1}. Suppose that
\begin{enumerate}
\item[(1)] There is a constant $\varepsilon>0$ such that $f\geq \varepsilon(\sum_{i=1}^n z_i\bar z_i)$ as $|z|\to \infty$.
\item[(2)] $\omega_f$ is complete.
\end{enumerate}
Then   the second derivatives of $f$ of mixed type are  constants.
\end{theorem}
We have following corollary.
\begin{corollary}\label{corollary_1.2} Let $f\in \mc R^\infty(\Omega)$ satisfying \eqref{eqn_1.1}. Suppose that
\begin{enumerate}
\item[(1)] There is a constant $C>0$ such that
$$  C^{-1}(\sum_{i=1}^n z_i\bar z_i)\leq f\leq C(\sum_{i=1}^n z_i\bar z_i),\;\;\;\;as\;|z|\to \infty$$
\item[(2)] $\omega_f$ is complete.
\end{enumerate}
Then    $f$ must be affine equivalent to a quadratical polynomial.
\end{corollary}

To state our second theorem we define a invariant
$$D:=\sum f^{i\bar k}f^{j \bar l}f_{,ij}f_{,\bar b\bar l},$$
where $``,"$ denote the covariant derivative with respect to the metric $\omega_{f}$.
\begin{theorem}\label{theorem_1.2} Let $f\in \mc R^\infty(\Omega)$ satisfying \eqref{eqn_1.1}. Suppose that
\begin{enumerate}
\item[(1)] There is a constant $C>0$ such that
$$  C^{-1}(\sum_{i=1}^n z_i\bar z_i)\leq f\leq C(\sum_{i=1}^n z_i\bar z_i),\;\;\;\;as\;|z|\to \infty$$
\item[(2)] $f$ defined on whole $\mathbb C^n$.
\item[(3)] There is a constant $C>0$ such that
$D\leq C(1+f)^2.$
\end{enumerate}
Then $f$ must be affine equivalent to a quadratical polynomial.
\end{theorem}
By the assumption $f\geq \varepsilon(\sum_{i=1}^n z_i\bar z_i)$ we have $f$ attain its minimal at some point $p^\bullet$.
  Consider the coordinate transformation
$\tilde z=z-z(p^\bullet)$. By adding a constant $f$ attains its minimum $\min f=0$ at $\tilde z(p^\bullet)=0.$ Then there exists $R_{1}>R_{0}$ such that $f\geq \frac{\varepsilon}{2}(\sum_{i=1}^n \tilde z_i\bar {\tilde z}_i)$ when $|\tilde z|\geq R_{1}$.
Without loss of generality we assume that $f$ attains its minimum $\min f$ at $0$, and $f(0)=0$.

\section{\bf Proof of Theorem \ref{theorem_1.1}}
\v

\v
\subsection{\bf Estimate of $T$} Denote $T=\sum_{i=1}^n f^{i\bar j}g_{i\bar j}$, where $g_{i\bar j}$ is the Euclidean metric on $\mathbb C^n$.
Denote by $R_{i\bar j}(\omega_{f}) $ and $R_{i\bar j}(\omega_{g}) $   the Ricci curvature of the metric $\omega_{f}$ and $\omega_{g}$ respectively.
Let $C>0$ be an arbitrary large number. Consider the geodesic ball $B_f(0, 2C)$ centered at $0$ with redius $r=2C$. 
We take a linear transformation
\begin{equation}\label{eqn_2.1}
\dot{f}=\frac{f}{C^2},\;\;\;\dot{g}=\frac{g}{C^2},\;\;,\;\;\dot{z}=\frac{z}{C}.\end{equation}
Note that
\begin{equation}\label{eqn_res_f}
 \frac{\p^2 f}{\p z_{i}\p \bar z_{j}} = \frac{\p^2 \dot f}{\p \dot z_{i}\p \bar {\dot z}_{j}},\;\;\; \dot{T}:=\sum \dot{f}^{i\bar j}\dot{g}_{i\bar j}=T.
\end{equation}
Denote  $ \dot{f}_{i\bar j}=\frac{\p^2 \dot f}{\p \dot z_{i}\p \bar {\dot z}_{j}}.$ Obviously,
$$
\det(\dot{f}_{i\bar j})=1.
$$

Denote    $\|\cdot \|_{  f}$,  $\|\cdot \|_{g}$ and  $\|\cdot \|_{\dot f}$,   the norm with respect to the metric $\omega_{f}$, $\omega_{g}$ and $\omega_{\dot f}$ respectively.
Denote by $\dot{B}_{\dot{f}}(0,2)$ the geodesic ball centered at $0$ with radius $2$ with respect to the metric $\omega_{\dot f}$.
We estimate
$\|\nabla \dot{f}\|^2_{\dot{f}}$. In \cite{CLS1} Chen, Li and Sheng proved
\begin{lemma}\label{lemma_2.1}
Let $\dot f\in \mc R^\infty(\Omega) $  with
$\dot f(0)=\inf_{\Omega} \dot f=0$. Suppose that
 \begin{equation}\label{eqn_2.2}
\;\; \|R_{i\bar j}(\omega_{\dot f})\|_{\dot f}\leq \mff N_1,
 \;\;in \;\; B_{\dot{f}}(0,2),
\end{equation}
where $R_{i\bar j}(\omega_{\dot f})$ is   the Ricci curvature of the metric $\omega_{\dot f}$.
Then in $  B_{\dot{f}}(0,1)$
\begin{equation}\label{eqn_2.3}
 \frac{\|\nabla \dot{f}\|_{\dot{f}}^2}{(1+ \dot{f})^2}\leq \mff C_{1}
\end{equation}
where  $\mff C_{1}>0$ is a constant depending only on $ n$ and
$N_0.$ Then,  for any $q\in B_{\dot{f}}(0,1),$
\begin{equation}\label{eqn 2.4}
  \dot{f}(q)-  \dot{f}(0) \leq \exp\left(\sqrt {\mff C}_{1}d(0,p)\right)-1,
   \end{equation}
   where $d(0,p)$ denotes the geodesic distance from $0$ to $p$ with respect to the metric $\omega_{\dot f}.$
\end{lemma}
\v
By the condition (1) of Theorem \ref{theorem_1.1} and the
linear transformation \eqref{eqn_2.1} we have,
 $$ \dot f\geq \varepsilon(\sum_{i=1}^n \dot z_i\bar{\dot z}_i).$$
 By Lemma \ref{lemma_2.1}, in $B_{\dot{f}}(0,1)$, we conclude that
\begin{equation}\label{eqn_2.5}
\sum |\dot{z}_i|^2\leq \frac{1}{\varepsilon}\exp(\sqrt {\mff C_{1}}).
\end{equation}

\v
To estimate $\dot T$, we need
\begin{lemma}\label{lemma_2.2}Let $\dot f  \in \mc R^\infty(\Omega)$
and $B_{\dot{f}}(0, 1)\subset \Omega.$  Suppose
\begin{equation*}det(\dot{f}_{i\bar j})\leq \mff N_1,\;\;\;  \|R_{i\bar j}(\omega_{\dot f})\|_{\dot f} \leq \mff N_1,
\;\;\; |\dot z|\leq
\mff N_1.\end{equation*}in $B_{\dot{f}}(0, 1)$, for some constant $\mff N_1>0$. Then
there exists a constant $\mff C_{2}>1$ such that
$$\mff C_{2}\inv\leq \lambda_1\leq \cdots
\leq\lambda_n\leq \mff C_{2},\;\;\forall\;q\in
B_{\dot{f}}(0, 1/2).$$
 where $\lambda_1,\cdots,\lambda_n$ are eigenvalues of the matrix $(\dot f_{i\b
 j}),$ $\mff C_{2}$
 is a positive constant depending on $n$ and $\mff N_1.$
  \end{lemma}
The Lemma \ref{lemma_2.2} has been proved in \cite{CLS1}.
By Lemma \ref{lemma_2.2} and $\dot g_{ i\bar j}=\delta_{i\bar j}$, $\dot{T}$ is bounded by a constant in
$B_{\dot{f}}(0,1/2)$.  By \eqref{eqn_res_f} we have  $T$ is bounded  by a constant in $B_{f}(0, C/2)$.
 Since $C$ is arbitrary,  $\sum f^{i\bar i}=\sum_{i} 1/\lambda_{i}$ and $\sum f_{i\bar i}=\sum_{i} \lambda_{i}$ we have
\begin{equation}\label{eqn_eqv_metric}
\mff C_{2}^{-1}\omega_{g}\leq \sum \omega_{f}\leq \mff C_{2}\omega_{g}.
\end{equation}
on whole K\"{a}hler  manifold $(\Omega, \omega_f)$.

\v
Set $\varphi=f-g.$  Denote
$$
S_{f} :=\sum f^{i\bar j}f^{k\bar l}f^{m\bar n}\varphi_{;i\bar l m}\varphi_{;\bar j k \bar n}.
$$
where $``;"$ denote the covariant derivative with respect to the metric $\omega_{g}$. If no danger of confusion we denote $S_{f}$  by $S.$
 It was proved that (see \cite{Y,PSS})
\begin{equation}\label{eqn_df_S}
\df S\geq -\mff C_{3}S-\mff C_{4}
\end{equation}
for some constants $\mff C_{3},\mff C_{4}>0,$
where $\df $ is the Laplacian operator with respect to the metric $\omega_{f}$.
\v
The following lemma can be found in \cite{CLS1}.
\begin{lemma}\label{lemma_5.2.2}
\begin{eqnarray} \label{eqn_df_T}
\df T&\geq &\sum
f^{i\bar{j}}f^{n\bar{l}}f^{k\bar{h}}f^{m\bar{p}}g_{k\bar{l}}
\phi_{;n\bar{p}\bar{j}}\phi_{;m\bar{h}i}  +\sum
f^{m\bar{l}}f^{k\bar{h}}g_{k\bar{l}}  R _{m\bar{h}}(\omega_{f}) \\&&- \|R_{i\bar j k\bar l}(\omega_{g})\|_{g} T^2,  \nonumber \\
 \label{eqn_df_logT}
 \df   \log T & \geq & -\| R_{i\bar j}(\omega_{f})\|_f - \|R_{i\bar j k\bar l}(\omega_{g})\|_{g}T.
\end{eqnarray}
where $R_{i\bar j k\bar l}(\omega_{g})$ denotes the  holomorphic bisectional curvature of the metric  $\omega_{g}$.
\end{lemma}

Since $g_{i\bar j}=\delta_{ij}$   and $R_{i\bar j}(\omega_{f})=0$, by \eqref{eqn_df_T} and \eqref{eqn_df_logT} we have
\begin{align}\label{eqn 2.4}
\df T\geq  \frac{1}{\lambda_{i}\lambda_{m}\lambda^2_{k}}|\phi_{;k\bar m i}|^2 \geq \mff C_{1}^{-1} S,\;\;\;\;\;\;\df \log\;T\geq 0.
\end{align}

\v
\subsection{\bf Estimate of $S$}
By \eqref{eqn_df_S} and \eqref{eqn 2.4}
we have
\begin{lemma}\label{lemma_2.3} There are constants $A>0$, $\mff C_{5},\mff C_{6}>0$ such that
$$\df (S+ AT)\geq  \mff C_{5}S -\mff C_{6}.$$
 \end{lemma}
Now we prove
\begin{lemma}\label{lemma_2.4} Suppose that $T\leq \mff C_{2}$ in $B_f(0,a)$, where $a\geq 1$. Then,
in $B_f(0,a/2)$, we have
$$S\leq \mff C_{7}\left(1 +\frac{1}{a^2}\right).
$$
where $\mff C_{7}$ is a constant depending only on $\mff C_{1}.$
\end{lemma}
\v\n
{\bf Proof.} Consider the function
$$\mathbf{F}=(a^2-r^2)^2e^{\alpha T}(S+ AT)$$
defined on $B_f(0,a)$, where $\alpha<1$ is a constant to be determined later. $\mathbf{F}$ attains its supremum at some interior point $p^\ast$.
Then, at $p^\ast$,
\begin{equation}\label{eqn_2.6}
\left(-\frac{2(r^2)_{,i}}{a^2-r^2}+\alpha  T_{,i}\right)
+\frac{ (S+ AT)_{,i}}{(S+AT)}=0,\end{equation}
\begin{equation}\label{eqn_2.7}
\left(-\frac{2\df (r^2)}{a^2-r^2} - \frac{2\|\nabla(r^2)\|_{f}^2}{(a^2-r^2)^2}+
\alpha \df T\right)
+\frac{ \df (S+ AT)}{(S+AT)}-\frac{ \|\nabla(S+ AT)\|_{f}^2}{(S+AT)^2}\leq 0. \end{equation}
Using the Schwarz inequality and \eqref{eqn_2.6} we obtain
\begin{equation}\label{eqn_2.8}
\frac{ \|\nabla(S+ AT)\|_{f}^2}{(S+AT)^2}\leq (1+\frac{1}{\delta})\frac{4 \|\nabla r^2\|^2}{(a^2-r^2)^2}+\alpha^2(1+\delta)  \|\nabla T\|_{f}^2
\end{equation}
On the other hand, it is easy to check that
 \begin{equation}\label{eqn_2.9}
\|\nabla T\|_{f}^2\leq n\mff C_{1}^2S.
\end{equation}
In fact, we choose local coordinates such that
$$
f_{i\bar j}=\lambda_{i}\delta_{ij},\;\;\;\;g_{i\bar j}=\delta_{ij}.
$$
Then
$$
\|\nabla T\|_{f}^2=\sum_{k}\frac{1}{\lambda _{k}} \left|\sum_{i}\frac{1}{\lambda_{i}^2}\varphi_{;i\bar i k}\right|^2 \leq \mff C_{1}^2 \sum_{k}\frac{1}{\lambda _{k}} \left|\sum_{i}\frac{1}{\lambda_{i}}\varphi_{;i\bar i k}\right|^2\leq n\mff C_{1}^2 S.
$$

Note that $|\nabla r|=1$. Using the Laplacian
comparison theorem we have
\begin{equation}\label{laolacian comparison}
\frac{2\df (r^2)}{a^2-r^2} + \frac{2\|\nabla(r^2)\|^2}{(a^2-r^2)^2}\leq \frac{4a^2}{(a^2-r^2)^2}+\frac{C(n)}{a^2-r^2}.\end{equation}
Then we obtain that
\begin{equation}
(1-n \mff C^{3}_{1}(1+\delta)\alpha)\mff C_{1}^{-1}\alpha S\leq C\left(1 +  \frac{a^2}{(a^2-r^2)^2} \right)
\end{equation}
Choose $\delta=1$ and $\alpha$ small such that
$n\mff C^{3}_{1}(1+\delta)\alpha\leq \frac{1}{2}.$ By $T<n\mff C_{1}$ we have
$$(a^2-r^2)^2e^{\alpha T}(S+AT)\leq C\left(a^{4} +   a^2\right).$$
Then Lemma \ref{lemma_2.4} is proved.
\v\v

\subsection{\bf Proof of Theorem \ref{theorem_1.1}} Let $p=(z_o, f(z_o))$ be an arbitrary point in $(\Omega, \omega_f)$, we prove $S_{f}(p)=0$. Suppose that $S_{f}(p)\ne 0$, we take a sequence $C_k$ and a sequence of geodesic ball $B_f(0, 2C_k)$ with $C_k\rightarrow \infty$. We take a sequence of linear transformations:
\begin{equation}
\dot{f}_k=\frac{f}{C^2_k},\;\;\dot{z}_k=\frac{z}{C_{k}}.\end{equation}
We get a sequence of geodesic ball $B_{\dot{f}_k}(0, 2)$.  Obviously, $p\in B_{\dot{f}_k}(0, 1/2)$ as $k$ large enough. Since $T$ is invariant under these   coordinate transformations, we have
$$S_{\dot{f}_k}(p)=C_{k}^2 S_{f}(p)\rightarrow \infty\;\;\;if\;\;S_{f}(p)\ne 0.$$
On the other hand, by $p\in B_{\dot{f}_k}(0, 1/2)$ and
  Lemma \ref{lemma_2.4}, $S_{\dot f_{k}}(p)$ are uniformly bounded above by a constant independent of $k$, we get a contradiction.
So $S_{f}\equiv 0$ on whole K\"{a}hler  manifold $(\Omega, \omega_f)$.
Theorem \ref{theorem_1.1} is proved.    \;\;\;$\blacksquare$
\v\v

\begin{remark}
By  \eqref{eqn_eqv_metric},  and the theorem of  Riebesehl and Schulz (cf. Theorem 2 in \cite{RS}) we can immediately obtain  Theorem \ref{theorem_1.1}.  Here we give a different proof, which has independent interest for us.
\end{remark}

\n{\bf Proof of Corollary \ref{corollary_1.2}.} By Theorem  \ref{theorem_1.1} we conclude that
$ f_{i\bar j},1\leq i,j\leq n$ are constants.  Set $v=\sum_{i,j}z_{i}f_{i\bar j}(0)\bar z_{j}.$ Then $f-v$ satisfying
$$
\frac{\p^2}{\p z_{i}\p \bar {z}_{j}}(f-v)=0,\;\;\;\;\;\forall 1\leq i,j\leq n,$$
$$    |f-v|\leq C (1+|z|^2).
$$
In particular, $f-v$ is a harmonic function.   By the  estimates of harmonic function (cf. Theorem 1.13 in \cite{HL})
we have, for any $R>1$ and  any multi-index $\nu$ with $|\nu|=2,$
$$|\nabla ^\nu (f-v)(p)|\leq  \frac{C(n)}{R^2}\max_{D(p,R)}|f-v|\leq \frac{C(n)}{R^2}\max_{D(0,R+d_{E}(0,p))}|f-v| ,$$
where $d_{E}(0,p)$ denotes the Euclidean distance from $0$ to $p,$ $D(p,R)$ denotes Euclidean ball centered at $p$ with radius $R$.
Choose $R$ big enough, we have $|\nabla ^\nu (f-v)(p)|\leq C.$
 Then by  Liouville Theorem we have $f-v$ is quadratic.    \;\;\;$\blacksquare$

\section{\bf Proof of Theorem \ref{theorem_1.2}}

Denote $S_{f}(0,b):=\{z | f(z)\leq b\}$. Consider the function defined on the section $S_{f}(0,b)$
$$\mathbf{F}:=\exp\left\{-\frac{b}{b-f}\right\}\frac{\|\nabla f\|_{f}^2}{(1+f)^2}.$$
$\mathbf{F}$ attains its supremum at some interior point $p^\ast$. We may assume $\|\nabla f\|_{f}(p^\ast)>0.$
Then, at $p^\ast$,
\begin{equation} \label{eqn_3.1}
-\left[\frac{2f_{,i}}{1+f}+\frac{b f_{,i}}{(b-f)^2}\right]\sum
f^{i\b j}f_{i}f_{\b j} +\sum f^{i\b j}
f_{,i  k}f_{\b j}+f_k=0.\end{equation}
 Choose a complex  coordinate system such
that, at $p^\ast$,
\begin{equation*}
f_{i\bar j}=\delta_{ij},\;\;f_{,1} = f_{,\bar{1}},\;\;\;f_{,i} =
f_{,\bar{i}}=0\;\;\forall \;i>1.\end{equation*}
Then from \eqref{eqn_3.1} we have
\begin{equation} \label{eqn_3.2}
-\left[\frac{2f_{,1}}{1+f}+ \frac{b f_{,1}}{(b-f)^2}\right]f_{,1}f_{,\bar 1} +
f_{,11}f_{\b 1}+f_1=0.\end{equation}
From \eqref{eqn_3.2} and the condition (3) of Theorem \ref{theorem_1.2} it follows that
$$\frac{2\|\nabla f\|_{f}^4}{(1+f)^2}\leq C(1+f)^2 + 1.$$
Then
\begin{equation} \label{eqn_3.3}
\exp\left\{-\frac{b}{b-f}\right\}\frac{\|\nabla f\|_{f}^2}{(1+f)^2}\leq \frac{\sqrt{C + 1}}{\sqrt{2}}.\end{equation}
Since $\mathbf{F}$ attains its supremum at $p^*$, \eqref{eqn_3.3} holds everywhere on $S_{f}(0,b)$.
Letting $b\to
\infty$, we have
\begin{equation} \label{eqn_3.4}
\frac{\|\nabla f\|^2}{(1+f)^2}\leq e\frac{\sqrt{C + 1}}{\sqrt{2}}:=\mff C_{9}.\end{equation}
Take a point  $p_1\in \p S_{f}(0,b)$ such that $d(0,p_1)=d(0,\p
S_{f}(0,b)).$ Let $l$ be the shortest geodesic from $0$ to
$p_1.$ Following from \eqref{eqn_3.4} we have $$\sqrt{\mff C_{9}}  \geq \|\nabla \log
(1+ f)\|  \geq  \left|\frac{d\log (1+f)}{ds}\right| ,$$ where $s$ denotes the arc-length parameter with
respect to the metric $f_{i\bar j}$. Applying this and a direct
integration we obtain
\[d(0,p_1)= \int_l ds \geq \mff C_{9}^{-\frac{1}{2}} \log(1+f)|_{p_1}=
 \mff C_{9}^{-\frac{1}{2}} \log\left(1+b\right) . \]
As $b\to \infty$, we obtain $d(0,\p S_f(0,b))\to +\infty.$
Hence
 $  \omega_f$ is complete. Then we use Corollary \ref{corollary_1.2} to complete the proof of Theorem \ref{theorem_1.2}.  $\blacksquare$

\v

{\bf Acknowledgment.}  We would like to thank Max-Planck-Institut
f\"ur Mathematik in den Naturwissenschaften, especially Professor J\"urgen Jost and Professor Xianqing Li-Jost,  for their great hospitality.

\v\v

\bibliography{<your-bib-database>}

\end{document}